# Multiscale simulation of injection-induced fracture slip and wing-crack propagation in poroelastic media

Hau Trung Dang*, Inga Berre, Eirik Keilegavlen

Center for Modeling of Coupled Subsurface Dynamics, Department of Mathematics, University of Bergen

*Corresponding author: Department of Mathematics, University of Bergen, Postboks 7803, 5020 Bergen, Norway

E-mail addresses: Hau.Dang@uib.no; dtrhau@gmail.com (H. Dang-Trung)

Inga.Berre@uib.no; Eirik.Keilegavlen@uib.no

**Abstract**

In fractured poroelastic media under high differential stress, the shearing of pre-existing fractures and faults and the corresponding propagation of wing cracks can be induced by fluid injection. This paper presents a two-dimensional mathematical model and a numerical solution approach for coupling fluid flow with fracture shearing and propagation under hydraulic stimulation by fluid injection. Numerical challenges are related to the strong coupling between hydraulic and mechanical processes, the material discontinuity the fractures represent in the medium, the wide range of spatial scales involved, and the strong effect that fracture deformation and propagation have on the physical processes. The solution approach is based on a multiscale strategy. In the macroscale model, flow in and poroelastic deformation of the matrix are coupled with the flow in the fractures and fracture contact mechanics, allowing fractures to frictionally slide. Fracture propagation is handled at the microscale, where the maximum tangential stress criterion triggers the propagation of fractures, and Paris' law governs the fracture growth processes. Simulations show how the shearing of a fracture due to fluid injection is linked to fracture propagation, including cases with hydraulically and mechanically interacting fractures.

## 1. Introduction

In the hydraulic stimulation of geothermal reservoirs in igneous rocks, elevated pressures in combination with anisotropic stress conditions result in shear displacement and the dilation of fractures and faults favorably oriented to slip, propagation of wing cracks from sliding or shearing fractures, and/or propagation of hydraulic fractures [1–4]. Sliding, dilation, and the propagation of fractures affect the stress and flow regime in the formation and, thereby, the stress state and deformation of nearby fractures. The coupling between flow in fractured and faulted rocks, fracture slip and propagation, and poromechanical matrix deformation is characterized by its multiscale nature: fracture propagation occurs locally but impacts and interacts with macroscopic reservoir-scale flow and deformation of the fractured rock.

In the current work, we present a modeling approach for hydraulic stimulation of fractured reservoirs under anisotropic stress. In this case, depending on the elevation of fluid pressure, the stimulation will cause slip of pre-existing fractures as well as fracture propagation. Slips of pre-existing fractures occur when coupled hydromechanical processes induced by fluid injection result in changes to the effective stress regime so that the fracture's frictional resistance to slip is exceeded [5,6]. The stress alterations resulting from fracture slip are coupled with fluid pressurization and drive tensile propagation of wing cracks at the fracture’s tips. Hence, in contrast to most of the research literature on fracture propagation resulting from hydraulic stimulation, we do not only consider the development of tensile hydrofractures. Instead, the reservoir stimulation is caused by a combination of slip of pre-existing fractures with fracture propagation [1–3]. Following McClure and Horne [2], we refer to this as mixed-mechanism stimulation.

To fully represent how injection operations alter fractured rock characteristics, simulation models must capture both the slip and deformation of existing fractures as well as fracture propagation. In addition, they must be able to account for the heterogeneous characteristics of subsurface formations. Challenges are related to capturing how the hydromechanical processes in the matrix interact with the flow, deformation, and propagation of fractures. This includes accounting for fracture contact mechanics, with the possibility of fractures being closed, sliding, and open.

To model the physics of these phenomena, fractures must be represented explicitly in an otherwise intact porous medium, conceptually leading to a discrete fracture-matrix model. To avoid resolving thin fractures in their normal direction, fractures are represented as co-dimension one objects embedded in the host medium with corresponding dimensionally reduced equations, resulting in a mixed-dimensional model [7]. Discretizations of such models can be both nonconforming and conforming. While conforming methods align the computational grid to the fractures, nonconforming methods utilize enrichment functions to capture the effects of the fractures. Recently, the nonconforming XFEM methodology, which has previously been developed for hydraulic fracture propagation [8–11], has been extended to couple fracture propagation with fracture contact mechanics in impermeable media [12,13]. A conforming method allows the direct assignment of variables and governing equations to the host medium, the fracture, and the matrix-fracture interface. For fluid flow, this flexibility simplifies the task of correctly capturing fluid exchange within and between fractures and the matrix. Furthermore, a conforming fracture representation allows for the modeling of fracture contact mechanics in poroelastic media [9,10,14–16], including fracture propagation [17,18] in a manner that directly couples shear and normal displacements along the fracture with alterations in stress and flow regime [19].

In the numerical modeling of the fractures mechanics in porous media caused by fluid injection, numerical models have generally focused either (1) on the tensile propagation of hydraulic fractures or (2) the deformation of pre-existing fractures or faults. The first regime ignores the deformation of existing fractures caused by frictional contact and governs the expansion of fracture by Mode I fracturing [8,10,15,16,20–31]. The second regime focuses on the frictional sliding, dilation and/or opening of pre-existing fractures or faults without considering fracture propagation [11,19,32–37]. Some recent models couple both regimes [2,18,19,38]. However, the coupling of hydromechanical processes between fractures and the matrix and the fracture propagation in these models are typically based on strong assumptions. For example, fracture dilation is assumed to not affect the surrounding stress [2,17], and restrictions allow fractures to propagate only along predefined paths [2] or the edges of fixed grid cells [18,38].

In this work, we present a numerical model, based on conforming discretizations, for injection-induced fracture shear-deformation and wing-crack propagation in poroelastic media without these limitations. Biot's model for poroelasticity governs flow and deformation in the matrix, with the deformation of fractures represented by contact mechanics, which consider the fractures' frictional resistance to sliding [19]. The fracture growth process is governed by the maximum tangential stress criterion and Paris's law [39,40]. A relatively coarse grid can be accepted for flow, poroelastic matrix deformation, and fracture deformation. In contrast, to capture the stress and correctly evaluate fracture propagation, a refined grid is needed around the fracture tip. As fracture propagation occurs locally from the fracture’s tips [2,41], an efficient solution strategy can be defined based on multilevel methods [42–44], in which the heterogeneous multiscale approach [42,45] is used in this work. Flow, poroelastic deformation, and contact mechanics of fractures are evaluated in a macroscale model, and mechanical fracture propagation is evaluated in a local microscale model [46] subject to body forces and boundary conditions that also account for the influence of macroscale fluid pressure. The models are coupled via displacement fields close to fracture tips (macroscale to microscale) and updates to the fracture path (microscale to macroscale model).

In discretizing the model on the coarse scale, a finite volume method for fracture and matrix flow and poroelastic matrix deformation is combined with an active set strategy for fracture contact mechanics [19]. For fracture propagation on a fine scale, a finite element method is applied in combination with collapsed quarter-point elements at the fracture tips to capture their stress singularity [17,47]. Adaptive remeshing is introduced on both scales to account for fracture propagation based on the implementation by Dang-Trung et al. [17]. We present numerical examples that focus on the method's ability to balance accuracy and computational cost under variations in grid resolution and the coupling between macroscale and microscale modeling. We also present a case with multiple hydromechanically coupled fractures, showcasing the capacity of our methodology to solve complex problems.

The paper is structured as follows. Section 2 presents the governing equations and Section 3 presents the multiscale solution strategy. Section 4 presents a numerical approach that employs a novel combination of a finite volume method for the poroelastic deformation of existing fractures with a finite element approach for the fracture propagation process. Section 5 presents several numerical test cases to show the stability and accuracy of the proposed approach and the potential in settings where multiple fractures mutually affect each other.

## 2. Mathematical model

The mathematical model for injection-induced fracture shear deformation and wing-crack propagation in poroelastic media is based on systems of partial differential equations and KKT conditions governing the physics. We start by introducing notations of geometry and primary variables. Then, we present the mathematical model for fracture contact mechanics, poroelastic deformation of the matrix, and fluid flow. At the end of this section, we present the model for the fracture propagation process.

### 2.1 Geometry and primary variables

As shown in Fig. 1, we represent a fractured porous media as a two-dimensional domain $\Omega$ that is divided into a host medium, termed the matrix and denoted $\Omega^{\mathrm{M}}$, and a set of fractures that are considered to be one-dimensional objects embedded in $\Omega^{\mathrm{M}}$. From here on, we will refer to both fractures and faults simply as fractures and let $\Omega_i^{\mathrm{F}}$ denote fracture $i$ and $\partial_k\Omega_i^{\mathrm{F}}, k = \{1, 2\}$ represent the two tips of $\Omega_i^{\mathrm{F}}$. Throughout the paper, we assume that fractures do not intersect. Finally, we denote the interface between $\Omega^{\mathrm{M}}$ and $\Omega_i^{\mathrm{F}}$ by $\Gamma_i$, where, when needed, we shall represent the two sides of the interface by $\Gamma_i^+$ and $\Gamma_i^-$. The boundary of $\Omega^{\mathrm{M}}$ that coincides with $\Gamma_i^{\pm}$ is denoted $\partial_i^{\pm}\Omega^{\mathrm{M}}$.

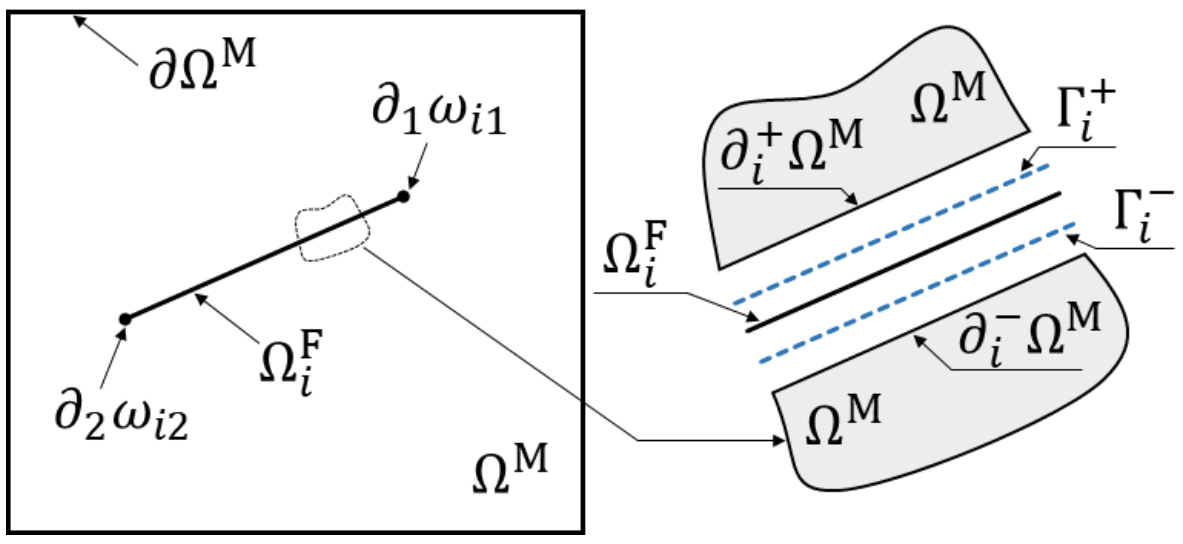

Fig. 1. Illustration of a fracture, $\Omega_i^F$ and a host medium, $\Omega^M$.

The primary variables are displacements, fluid pressures, contact forces on the fractures, and fluid fluxes between fractures and the matrix. Specifically, the displacements on $\Omega^{\mathrm{M}}$ are denoted as $\mathbf{u}$, the pressure is represented by $p$ in $\Omega^{\mathrm{M}}$ and $p_i$ in $\Omega_i^{\mathrm{F}}$, and the contact force $\mathbf{f}_i$ is defined only in $\Omega_i^{\mathrm{F}}$. Finally, $\mathbf{u}_i$ and $\lambda_i$ denote the displacement and fluid flux on $\Gamma_i$, respectively, where $\mathbf{u}_i$ can be split into $\mathbf{u}_i^+$ and $\mathbf{u}_i^-$ on $\Gamma_i^+$ and $\Gamma_i^-$, respectively, and $\lambda_j$ can be similarly split. Time derivatives are denoted by a dot, e.g., $\dot{\mathbf{u}}_i$.

## 2.2 Fracture contact mechanics

Let $\mathbf{n}_i$ be the normal vector to the fracture surface $\Omega_i^{\mathrm{F}}$, pointing from $\Gamma_i^+$ to $\Gamma_i^-$. We define the jump operator acting on $\mathbf{u}_i$ by

$$[\![\mathbf{u}_i]\!] = ([\![\mathbf{u}_i]\!]_{\mathrm{n}}, [\![\mathbf{u}_i]\!]_{\tau}) = \mathbf{u}_i^- - \mathbf{u}_i^+, \tag{1}$$

where $[\![\mathbf{u}_i]\!]_{\mathrm{n}}$ and $[\![\mathbf{u}_i]\!]_{\tau}$ denote the normal and tangential components of the displacement jump, respectively. In the normal direction, we require nonpenetration and define the normal component of the contact traction, $\mathbf{f}_{i,\mathrm{n}}$, to be negative in contact. This assumption gives rise to the KKT condition, i.e.

$$[\![\mathbf{u}_i]\!]_{\mathrm{n}} - g \geq 0, \quad \mathbf{f}_{i,\mathrm{n}} \leq 0, \quad ([\![\mathbf{u}_i]\!]_{\mathrm{n}} - g)\mathbf{f}_{i,\mathrm{n}} = 0. \tag{2}$$

Here, the gap function $g$ allows the fracture to open while the walls are still in mechanical contact. We set $g = [\![\mathbf{u}_i]\!]_\tau \tan\psi$, with $\psi$ being the dilation angle to let the fracture open due to tangential slip.

The tangential motion of the fracture is modeled as a frictional contact problem, with the following relation between the tangential contact traction, $\mathbf{f}_{i,\tau}$, and the change of displacement jump in time, $[\![\dot{\mathbf{u}}_i]\!]_\tau$:

$$
\begin{aligned}
&|\mathbf{f}_{i,\tau}| \leq -\mu_s \mathbf{f}_{i,\mathrm{n}}, \\
&|\mathbf{f}_{i,\tau}| < -\mu_s \mathbf{f}_{i,\mathrm{n}} \quad \text{so} \quad [\![\dot{\mathbf{u}}_i]\!]_\tau = 0, \\
&|\mathbf{f}_{i,\tau}| = -\mu_s \mathbf{f}_{i,\mathrm{n}} \quad \text{so} \quad \exists \gamma \in \mathbb{R},\ \mathbf{f}_{i,\tau} = -\gamma^2 [\![\dot{\mathbf{u}}_i]\!]_\tau.
\end{aligned} \tag{3}
$$

The tangential traction is bounded from above by the normal traction scaled by the friction $\mu_s$, and when the frictional resistance is overcome, the displacement is parallel to the tangential traction.

The force balance on the fracture's surfaces is given by

$$
\begin{aligned}
\mathbf{n}_i \cdot \boldsymbol{\sigma}|_{\partial_i^+ \Omega^{\mathrm{M}}} &= \mathbf{f}_i - \alpha_i p_i (\mathbf{I} \cdot \mathbf{n}_i), \\
-\mathbf{n}_i \cdot \boldsymbol{\sigma}|_{\partial_i^- \Omega^{\mathrm{M}}} &= \mathbf{f}_i - \alpha_i p_i (\mathbf{I} \cdot \mathbf{n}_i),
\end{aligned} \tag{4}
$$

where $\boldsymbol{\sigma}$ denotes the hydromechanical stress in the matrix, $\mathbf{f}_i$ is the contact traction acting on the fracture surface, and $p_i$ is the pressure inside the fracture. The Biot coefficient in the fracture is denoted by $\alpha_i$ and $\mathbf{I}$ is the identity matrix. Equality should be enforced on $\Gamma_i^{\pm}$, but for notational convenience, we have suppressed projection operators. See Keilegavlen et al. [48] for more information.

## 2.3 Flow and poroelastic deformation of matrix

Flow and deformation in $\Omega^{\mathrm{M}}$ are modeled by Biot theory, with the matrix taken as a linearly elastic medium. By neglecting inertial terms, the conservation of momentum and mass is governed by

$$
\nabla \cdot \boldsymbol{\sigma} = \mathbf{b}, \tag{5}
$$

$$
\alpha \frac{\partial (\nabla \cdot \mathbf{u})}{\partial t} + M \frac{\partial p}{\partial t} - \nabla \cdot \left( \frac{\mathcal{K}}{\mu} \nabla p \right) = q, \tag{6}
$$

where $\boldsymbol{\sigma}$ is the hydromechanical stress in $\Omega^{\mathrm{M}}$, defined by

$$
\boldsymbol{\sigma} = \mathbf{C} \nabla_s \mathbf{u} - \alpha p \mathbf{I}. \tag{7}
$$

Here, $\nabla_s$ represents the symmetrized gradient, $\mathbf{C}$ is the stiffness matrix, $\mathbf{b}$ denotes body forces, and $q$ is the fluid source term. The Biot coefficient of the matrix is $\alpha$, the Biot modulus is given by $M = \left( \phi c_p + \frac{\alpha - \phi}{K} \right)$, $c_p$ is the fluid compressibility, $\phi$ is the matrix porosity, $K$ is the bulk modulus, $\mathcal{K}$ denotes the permeability of the porous matrix, which is assumed to be isotropic, and $\mu$ is the fluid viscosity. On $\partial_i^{\pm} \Omega^{\mathrm{M}}$, continuity of the displacements is enforced so that $\operatorname{tr} \mathbf{u}|_{\partial_i^{\pm} \Omega^M} = \mathbf{u}_i^{\pm}$, where tr is the trace operator. As seen from $\Omega^{\mathrm{M}}$, the interface displacement thus acts as a Dirichlet boundary condition. The interface fluid flux $\lambda_i$ enters as a Neumann condition to the mass conservation equation.

### 2.4 Fluid flow in fractures and matrix-fracture interaction

By using the discrete fracture-matrix model, the fracture is represented explicitly in the domain. Following Stefansson et al. [19], who extended the work of Martin et al. [49] to discrete fracture-matrix models with changing apertures, the conservation of mass in fracture *i* is given by

$$\frac{\partial a_i}{\partial t} + a_i c_p \frac{\partial p_i}{\partial t} - \nabla \cdot \left( \frac{\mathcal{K}_i}{\mu} \nabla p_i \right) + (\lambda_i^+ + \lambda_i^-) = q_i. \tag{8}$$

Here, we assume that the fracture can be completely occupied by the fluid, i.e., that the fracture porosity and Biot's coefficient are equal to one. $\mathcal{K}_i$ is the fracture tangential transmissivity. The aperture $a_i = a_i^0 + [\![\mathbf{u}_i]\!]_n$ is computed by a sum of an initial value $a_i^0$ and an update due to fracture deformation. Therefore, the first term in Eq. (8) represents volume changes due to changes in aperture. The fracture transmissivity is related to aperture by the so-called cubic law, $\mathcal{K}_i = a_i^3/12$ [50]; i.e., $a_i$ equals the hydraulic aperture of the fracture. The term $(\lambda_i^+ + \lambda_i^-)$ represents inflow from the matrix over $\Gamma_i^\pm$, where $\lambda_i^\pm$ is the interface flux between the matrix and the fracture defined as

$$\lambda_i^\pm = -\kappa_i (p_i - \operatorname{tr} p^\pm), \tag{9}$$

where $\kappa_i = 2\mathcal{K}_i/(\mu a_i^2)$ is an expression of permeability normal to the fracture, $p^\pm$ represents pressures from the matrix at the two sides of the fracture, and it is understood that the fracture and matrix pressures should be projected onto $\Gamma_i^\pm$.

### 2.5 Fracture propagation

The propagation criterion is based on a criterion on maximum tangential stress [39] for mixed-mode fracturing. The theory postulates propagation when the maximum tangential stress in the process zone around a fracture tip exceeds a critical value. The direction of propagation is that of the maximum tangential stress. The tangential stress around a fracture tip can be expressed in polar coordinates as

$$\sigma_\theta^{\mathrm{L}}(r,\theta) = \frac{1}{\sqrt{2\pi r}}\left(K_{\mathrm{I}}\cos^3\frac{\theta}{2} - \frac{3}{2}K_{\mathrm{II}}\cos\frac{\theta}{2}\sin\theta\right), \tag{10}$$

and the crack grows in the direction $\theta_0$ if $\sigma_\theta^{\mathrm{L}}(r,\theta_0) = \frac{K_{\mathrm{IC.}}}{\sqrt{2\pi r}}$, where $K_{\mathrm{IC.}}$ is the fracture toughness. The propagation angle is given by

$$\theta_0 = 2\tan^{-1}\left(\frac{K_{\mathrm{I}}}{4K_{\mathrm{II}}} \pm \frac{1}{4}\sqrt{\left(\frac{K_{\mathrm{I}}}{K_{\mathrm{II}}}\right)^2 + 8}\right), \tag{11}$$

subjected to a condition

$$K_{\mathrm{II}}\left(\sin\frac{\theta_0}{2} + 9\sin\frac{3\theta_0}{2}\right) < K_{\mathrm{I}}\left(\cos\frac{\theta_0}{2} + 3\cos\frac{3\theta_0}{2}\right), \tag{12}$$

where $r$ is the distance from the tip. $K_{\mathrm{I}}$ and $K_{\mathrm{II}}$ are the stress intensity factors (SIFs). The propagation length can, in general, be computed by a Paris-type law [40], which, under the assumption that there is a single fracture inside every microscale domain, simplifies to a propagation length equal to a preset value, $l_{\max}$. If more than one crack grows simultaneously, then the tips with the higher energy in the fracture set advance further than the others. The increment for each tip is defined by

$$l_{\mathrm{adv}}^i = l_{\max}\left(\frac{G_i}{\max(G_i)}\right)^{0.35}, \tag{13}$$

where $G_i$ is the energy release rate for the $i^{\mathrm{th}}$ propagation crack [51].

## 3. Multiscale solution strategy

Our goal is to define a computational approach for the interaction between, on the one hand, deformation of fractures and domains due to hydromechanical stresses and, on the other hand, the propagation of fractures. While a fully coupled approach to the governing equations presented in the previous section is possible, it is impractical for three reasons. First, compared to the time scale of fluid flow, wing-crack propagation can be considered quasi-static and can therefore be loosely coupled to the fluid flow problem [52,53]. Second, there will be no wing-crack propagation for significant periods, and efficiency can be gained only by considering the macroscale problem. Third, while it is critical for the propagation problem to accurately capture the stress in the vicinity of the propagating fracture, including the singularity at the tip, the mechanical response of the more expansive reservoir, including other (potentially propagating) fractures, can be given a coarser representation.

Motivated by these observations, our computational model is based on a multiscale approach. We define subproblems for, on the one hand, large-scale fluid flow and deformation of fractures and the matrix, and on the other hand, fracture propagation and accompanying deformation locally around fracture tips. Reflecting the small time and length scales involved in propagation, we assign separate, microscale domains around each fracture tip to be used in the purely mechanical propagation calculation, referred to as the microscale problem. Conversely, the whole simulation domain is referred to as the macroscale domain, on which we solve the macro problem consisting of fluid flow and fracture and matrix deformation.

### 3.1 Macroscale and microscale models

The macroscale problem is defined on the geometry presented in Section 2, with the governing equations given by Eqs. (1) - (9), that is, frictional contact mechanics coupled with hydromechanical deformation in the matrix and on the fracture domains. Fracture propagation is not explicitly accounted for in the macroscale model but is instead updated from the solution to the microscale problem, as discussed in Section 3.2.

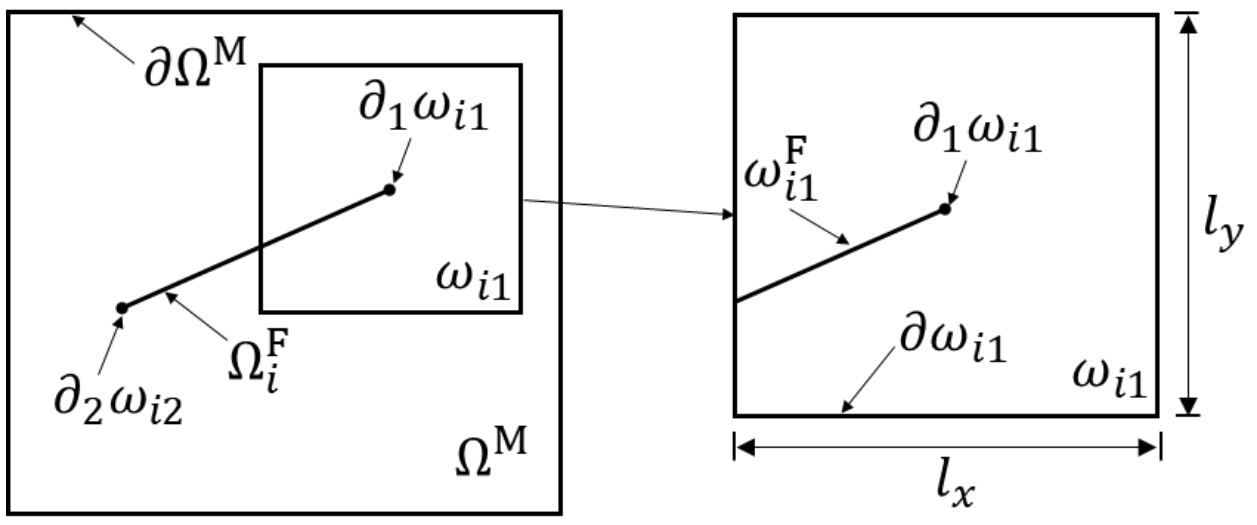

Fig. 2. Illustration of a fracture, $\Omega_i^F$ and a microdomain $\omega_{ik}$.

As shown in Fig. 2, the microscale models are centered on fracture tips. With each fracture tip $\partial_k\Omega_i^{\mathrm{F}}$ in the macroscale domain, we associate a (generally) smaller domain, termed a microdomain, and denoted as $\omega_{ik}$ with size $l_x \times l_y$. The microdomain is composed of a part of the matrix, $\omega_{ik} \subset \Omega^{\mathrm{M}}$, and a single fracture domain, $\omega_{ik}^{\mathrm{F}}$, which represents a part of the macroscale fracture $\Omega_i^{\mathrm{F}}$. In general, $\omega_{ik}^{\mathrm{F}} \not\subset \Omega_i^{\mathrm{F}}$ since the resolution of the propagating fracture is different in the micro and macro domains, as discussed in Section 4. We let $\partial_{\mathrm{F}}^{\pm}\omega_{ik}$ represent the two sides of the microscale fracture while $\partial_{\mathrm{M}}\omega$ is the rest of the boundary of $\omega$.

In accordance with the discussion at the beginning of this section, we include the effect of fluid pressure as a body force in the microscale model. The primary variable in the microscale problem is, therefore, the displacement in $\omega_{ik}$, which we represent by $\mathbf{u}^{\mathrm{L}}$ for simplicity. As with the full deformation, we assume that the microscale matrix behaves similarly to a linearly elastic and isotropic medium governed by

$$\nabla \cdot \mathbf{c}\nabla_s \mathbf{u}^{\mathrm{L}} = \mathbf{b} \tag{14}$$

where $\mathbf{b} = \nabla \cdot (\alpha p \mathbf{I})$ is the body force caused by pressure in the macroscale domain. $\mathbf{c}$ is the stiffness tensor. Boundary conditions for the microscale problem are set according to the macroscale state, as discussed next.

### 3.2 Coupling between macroscale and microscale models

For the microscale problem, the fracture surfaces are not allowed to move freely. Instead, their displacement is set from the macroscale behavior close to the fracture tip. That is, the microscale boundary condition on $\partial_{\mathrm{F}}^{\pm}\omega_{ik}$ is given as a displacement jump computed from the macroscale state,

$$[\![\mathbf{u}^{\mathrm{L}}]\!]|_{\partial_i^{\pm}\omega_{ik}} = \mathcal{R}\left([\![\mathbf{u}_i]\!]|_{\partial_i^{\pm}\Omega^{\mathrm{M}}}\right), \tag{15}$$

where $\mathcal{R}$ is a reconstruction operator, defined for the discrete problem in Section 4.4. On the remainder of the boundary, $\partial\omega_{ik}$, we similarly set

$$\mathbf{u}^{\mathrm{L}}|_{\partial\omega_{ik}} = \mathcal{R}(\mathbf{u}). \tag{16}$$

The coupling from the microscale to macroscale model consists of updating the macroscale fracture geometry based on microscale fracture propagation. This is linked to the continuous representation of the evolving geometry in the two models. For simplicity, we will represent the fractures as piecewise linear objects with a resolution related to that of the computational grids on the two scales, as detailed in Section 4. However, to cover more advanced features, including merging of fractures and three-dimensional problems, more elaborate geometric representations are needed [54,55].

## 4. Numerical approach

This section describes the building blocks of our numerical approach in terms of grids and solution approaches to macroscale and microscale problems, together with the multiscale coupling concept.

### 4.1 Numerical grids for fracture propagation

We construct numerical grids to conform to fractures in both the macroscale and microscale domains. That is, the grids on both $\Omega^{\mathrm{M}}$ and $\omega_{ik}$ are constructed so that fractures coincide with the paths of grid faces and then split nodes and faces along these paths, as done by Dang-Trung et al. [17]. In the macroscale domain, we further construct one-dimensional grids on $\Omega_i^{\mathrm{F}}$, as well as on the interfaces $\Gamma_i^{\pm}$.

Focusing on a single microscale domain, we represent the mesh size in a microscale grid by $\Delta h$ and let $\Delta H$ represent the macroscale mesh size around the same fracture tip. To ensure the stability of the propagation, the resolution of the microscale domain is set to be finer than that of the macroscale domain, i.e., $\Delta h = \varepsilon_{\mathrm{m}} \Delta H$ with $\varepsilon_{\mathrm{m}} \leq 1$. On both scales, we consider simplex grids, with the initial grids being constructed by Gmsh [56]. Updates to the fracture geometry will generally not follow existing paths of grid faces at both the micro- and macroscales. Thus, before fracture propagation on either the micro- or macroscale domain, the grid is adjusted in the vicinity region of radius $5 \times l_{\max}$ around the fracture tip with triangular rosette elements, size $\Delta h$ or $\Delta H$, to accommodate the extension of the fracture, followed by Laplacian smoothing to preserve the grid quality [17]. This technique is effective in removing degenerate and small elements. The fracture is then prolonged by splitting grid faces and nodes. For the macroscale grid, it is further necessary to prolong the grid for the fracture domain $\Omega_i^{\mathrm{F}}$ and the interfaces $\Gamma_i^{\pm}$ and to update the projections between the different grids.

We remark that although the grid adjustment necessitates an update of the discretization in $\Omega^{\mathrm{M}}$, the cost of this operation can be limited by confining the adjustment to a region close to the fracture tip.

### 4.2 Macroscale discretization

The governing equations (1) - (9) on the macroscale, namely poroelastic deformation in $\Omega^{\mathrm{M}}$ and both fluid flow and contact mechanics in $\Omega_i^{\mathrm{F}}$ and over $\Gamma_i$ are discretized and solved fully coupled by the open-source software tool PorePy [48]. The overall approach has been used before to study poroelastic [32,33] and thermoporoelastic [19] deformation coupled with fracture mechanics and has also been applied to field studies [57].

The discretization of the frictional contact problem requires handling the discontinuity in the contact conditions Eqs. (2) - (4). These are evaluated cellwise to determine whether fractures are open or closed and, if closed, whether they are sticking or slipping; see Stefansson et al. [19] for details. This classification is employed in an active set method, where the contact conditions and balance of forces expressed are discretized according to the state from the previous iteration [19,33,58].

The conservation equations for flow in $\Omega^{\mathrm{M}}$ and $\Omega_i^{\mathrm{F}}$ as well as momentum in $\Omega^{\mathrm{M}}$ are discretized by a family of cell-centered multipoint finite volume methods developed for poroelasticity [59,60]. The methods are based on constructing discrete representations of stresses (respective fluxes) over cell faces regarding displacements (respectively pressures) in nearby cell centers. The balance of momentum and mass is enforced on the cells. For fracture domains, the method reduces to the well-known two-point flux method, which can also deal with nonplanar domains resulting from fracture propagation. Finally, the coupling between the matrix and fractures follows the scheme described by Nordbotten et al. [14] for the flow problem.

The coupled set of equations is nonlinear and requests an iteration solver due to the active set approach to the contact conditions. The system is solved by a semismooth Newton method using a direct solver for the linearized system.

### 4.3 Microscale discretization

While our macroscale discretization was chosen to comply with conservation, calculations meant to decide whether a fracture will propagate, and if so, where it will go, pose different requirements on the spatial discretizations. Specifically, it is crucial to represent the stress singularity at the fracture tip. To that end, Eq. (14) is discretized by a finite element method with $\mathcal{P}_2$ basis functions. The stress singularity is captured using the nodal displacement correlation technique [61] based on quarter-point elements [47]. To enhance computational accuracy, the microscale grid is refined and guided by residual-based a posteriori error estimates [62]. From a computed displacement field, SIFs are estimated to determine whether the fracture should propagate and, if so, in which direction [17,61]. The microscale grid is updated as described in Section 4.1, and the displacement and pressure variables are mapped to the new grid by a $\mathcal{P}_1$ interpolation. For details of the algorithm and investigations of its performance on SIFs convergence and fracture propagation verification, we refer to Dang-Trung et al. [17].

Depending on the boundary conditions, several propagation steps may be needed to arrive at a stable state. During these iterations, the boundary conditions are fixed, consistent with the assumption that microscale propagation is instantaneous relative to dynamics on the macroscale.

### 4.4 Discrete mapping between macroscale and microscale models

To couple the numerical states on the macro- and microdomains, it is necessary to project displacements from the macro- to microscale domain boundaries and compress microscale updates to the fracture geometry onto the macroscale grid.

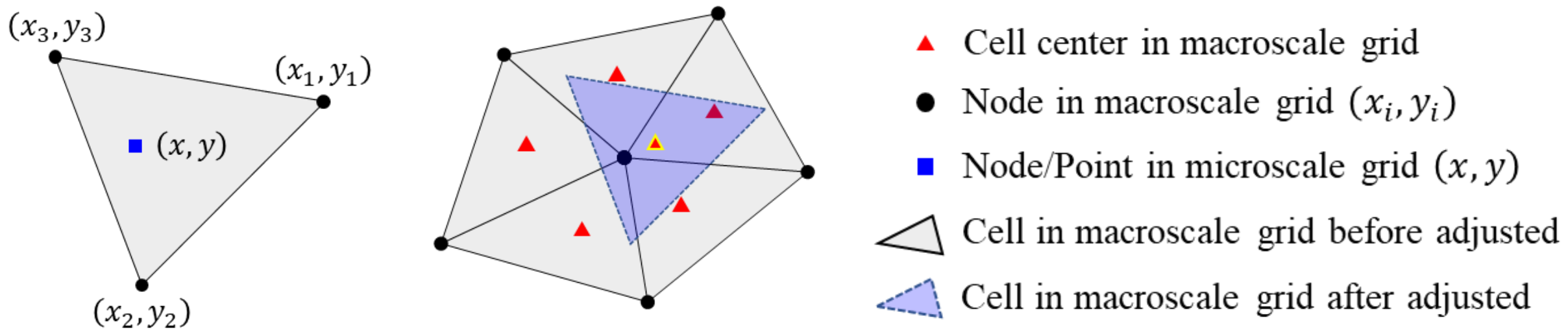

Fig. 3. Illustration of a point and grid in the micro and macroscale.

The values at any point in the microscale domain are determined through the values at the macroscale nodes by $\mathcal{P}_1$ interpolation. As illustrated in Fig. 3, for a microscale node with coordinates $(x, y)$ belongs to a parent cell defined by nodes $(x_i, y_i)$ in the macroscale grid, the values at $(x, y)$ is then approximated by

$$\xi(x,y) = \sum_{i=1}^{3} N_i(x,y)\xi(x_i,y_i), \tag{17}$$

where $N_i(x, y)$ are the Lagrange basis polynomials, defined by

$$\begin{aligned} N_1(x,y) &= \frac{1}{2A_e}[(y_2 - y_3)x + (x_3 - x_2)y + x_2y_3 - x_3y_2], \\ N_2(x,y) &= \frac{1}{2A_e}[(y_3 - y_1)x + (x_1 - x_3)y + x_3y_1 - x_1y_3], \end{aligned} \tag{18}$$

$$N_3(x,y) = \frac{1}{2A_e}[(y_1 - y_2)x + (x_2 - x_1)y + x_1y_2 - x_2y_1].$$

where $A_e$ is the area of the parent cell.

Using the FVM in the macroscale discretization, the displacements and pressure are naturally determined at the cell centers of the macroscale grid. So, we need to reconstruct node values from cell center values for calculation in Eq. (17). Besides, if a fracture propagates in the macroscale domain based on microscale solution, the macroscale grid is locally adjusted to ensure the fracture path coincides with faces and, therefore, the solution at cell centers on the adjusted grid also needs to be reconstructed. We deal with both cases by using the natural neighbor interpolation, such as

$$\xi_e = \frac{1}{\sum_{j=1}^{n} A_{ej}} \sum_{j=1}^{n} A_{ej}\xi_j, \tag{19}$$

where $\xi_e$ and $\xi_j$ denote values at the center of cell $e$ and $j$, and $A_{ej}$ is the intersection area.

Eq. (17) represents a discrete representation of the reconstruction operator $\mathcal{R}$ shown in Eqs. (15) and (16), i.e.,

$$\xi(x,y) = \tilde{\mathcal{R}}\xi_j^{\text{cell-center}} \tag{20}$$

where $\tilde{\mathcal{R}}$ is the discrete reconstruction. By Eq. (20), the values at desired microscale points are estimated through values at cells centers in the macroscale grid.

In our model, the fracture path is represented directly in the computational grid. Hence, the information transfer on the fracture path from the microscale to the macroscale model is dependent on the grid resolution on the two scales. Critical for simulation efficiency, small increases in the fracture length on the microscale domain are not immediately projected to the macroscale problem. Instead, for a fracture propagating on the microscale, with added length $|\Delta\omega_{ik}|$, the macroscale fracture is updated only when $|\Delta\omega_{ik}^{\mathrm{F}}| \geq \varepsilon_{\mathrm{p}}\Delta H$, where $\varepsilon_{\mathrm{p}}$ is a simulation parameter. When this threshold is overcome, the macroscale fracture is extended by a linear approximation of $\Delta\omega_{ik}^{\mathrm{F}}$ and the macroscale grid is updated as discussed in Section 4.1. Thus $\varepsilon_{\mathrm{p}}$ controls both the resolution of the macroscale grid in the vicinity of propagating fractures and the numerical coupling strength between the microscale and macroscale models.

### 4.5 Numerical solution approach

As a summary of the above presentation, Fig. 4 illustrates the workflow of the multiscale simulation approach. The time step size is usually taken as a constant represented by $\Delta t$. However, when a fracture propagates on the macroscale, both governing equations and parameters change along the fracture path, and the macroscale state adjusts accordingly. In particular, the fractured part of the rock experiences enhanced permeability and volume available for fluids. The pressure field in the vicinity of the tip will adjust to the new parameters on a time scale that is much shorter than that of pressure diffusion related to the injection. We capture this effect by temporarily reducing the time step size with a factor $\varepsilon = 10^{-2}$ until this rapid dynamics is resolved and the pressure field near the crack tip is stabilized, whereupon we continue with the standard step size.

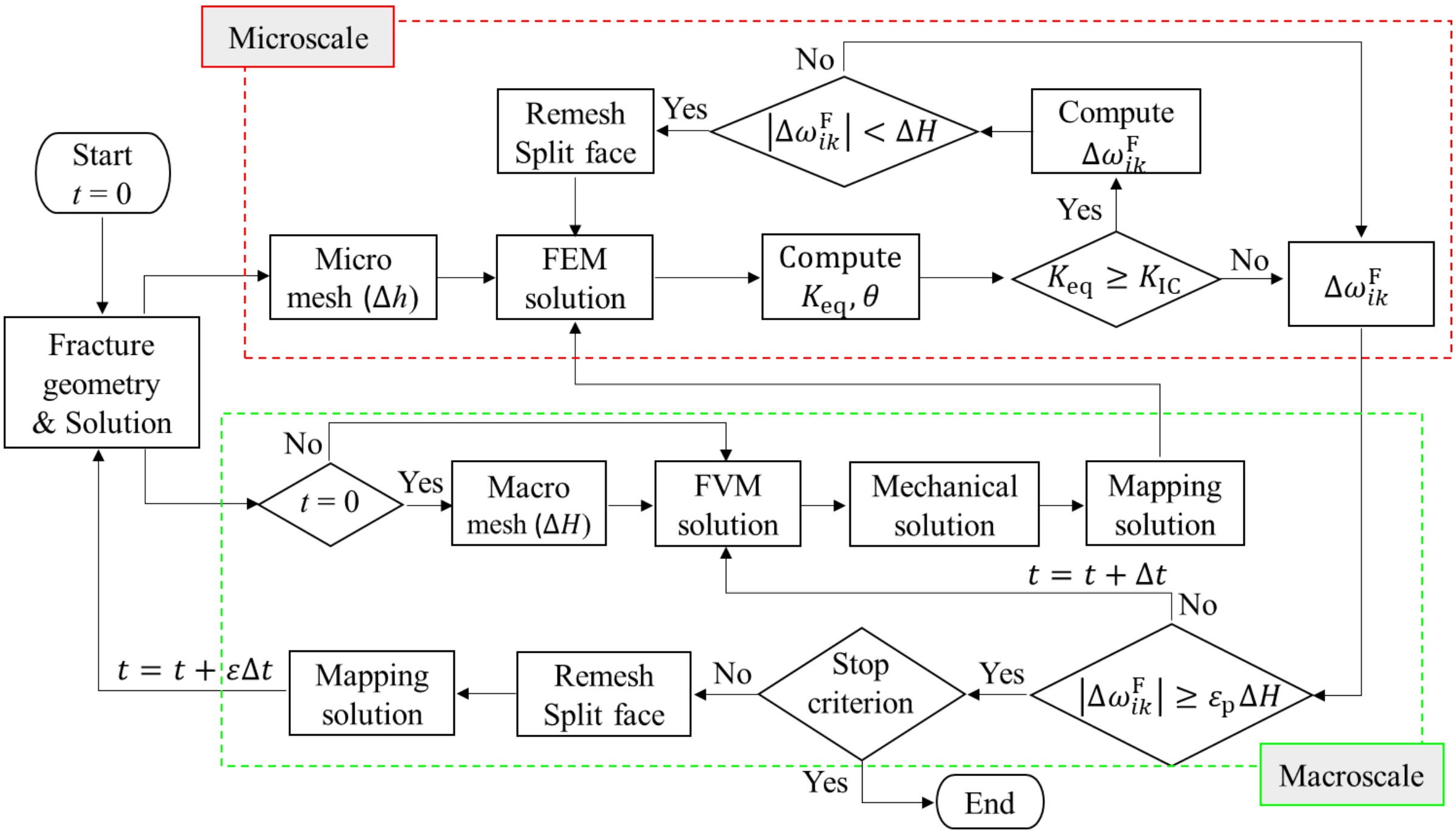

Fig. 4. Illustration of the workflow in the multiscale simulation method. The simulation is controlled by the macroscale mesh size $\Delta H$, time step size $\Delta t$, microscale domain size $l$, relation between microscale and macroscale grid size $\varepsilon_m$ and macroscale resolution of the fracture $\varepsilon_p$. The details of the "Remesh/Split face" box are given in Section 4.1.

## 5. Results

The correctness of either the micro or macro model has been verified in previous studies [5,17,19,33]. This section, therefore, is devoted to the presentation of numerical experiments of the fully coupled model. The macroscale problem alone may strain the available computational resources in application-oriented simulations with large domains and multiple fractures. Therefore, it is paramount to limit the additional computational cost to incorporate fracture propagation in such simulations. We have devised a suite of numerical experiments designed to investigate the stability, accuracy, and computational efficiency of the proposed numerical approach. Specifically, we study how the prediction of the fracture path is altered under variations in the size of the microscale domain ($l = l_x = l_y$), the mesh size on the microscale ($\Delta h$) and macroscale ($\Delta H$) domains, the time step size ($\Delta t$), and the threshold for updating the macroscale geometry ($\varepsilon_\mathrm{p}$). Together, these simulation parameters determine the balance between solution accuracy and computational cost. As our intention is to allow simulations on large domains where high-resolution simulations are not feasible, our focus is not on the convergence of the numerical solution but rather its stability as the resolution in time and space are coarsened. Based on observations from these tests, we finally present a complex case of hydromechanical processes interacting with the deformation and propagation of multiple pre-existing fractures in a synthetic subsurface fluid injection scenario. The source code for the following simulations can be found by referring to [63].

### 5.1 Onset of fracture

We first investigate how the onset of fracturing, determined by SIFs, is influenced by the size of the time step, macroscale mesh, and microscale domain. To that end, we consider a porous media domain assumed to be homogenous and linearly elastic with the material properties given in Table 1. The geometry and boundary conditions of the model are illustrated in Fig. 5, in which $L_x = L_y = 2$ m. A single fracture of length $l_f = 0.1$ m and initial aperture $a^0 = 1$ mm is located at the center of the computational domain and oriented at $45^\circ$ to the positive $x$-direction. The left and bottom boundaries are fixed in the $x$- and $y$-directions, respectively, and the top and right boundaries are free. The fluid is not allowed to flow through boundaries. The fractured porous media is subject to a stress state with the maximum horizontal stress $\sigma_1 = 20$ MPa and the minimum vertical stress $\sigma_2 = 10$ MPa, imposed orthogonally to the domain along the *x*- and *y*-directions, respectively. Water is injected into the pre-existing fracture after 6 hours (h) and continuously for 15 h at a constant rate of $Q_0 = 5 \times 10^{-9}\ \mathrm{m^2/s}$.

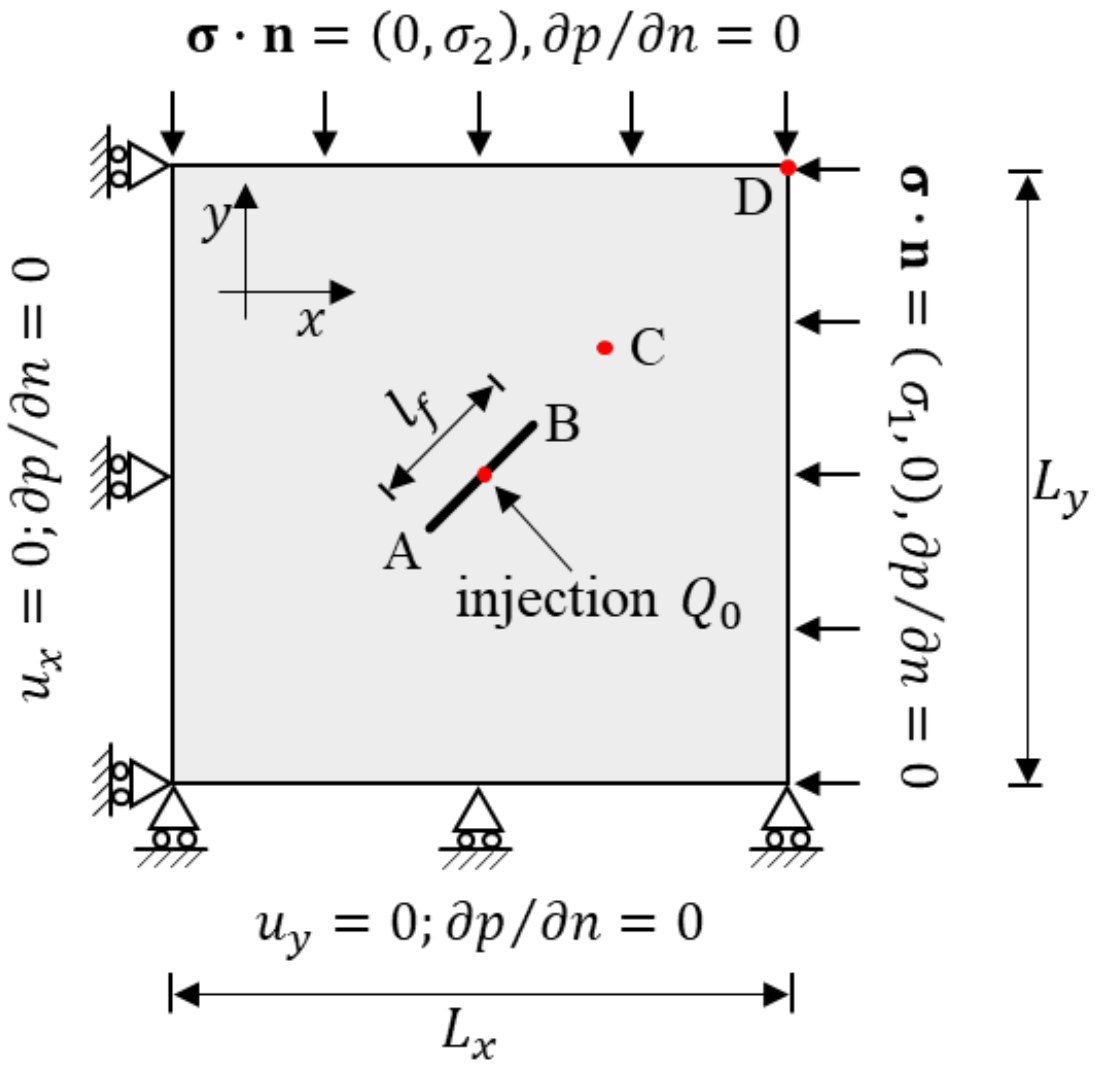

Fig. 5. Model of a porous media with a single fracture subjected to a principal stress regime.

Table 1. Material properties.

| Parameter | Definition | Value |
|---|---|---|
| $E$ | Young's modulus | 40 GPa |
| $\nu$ | Poisson's ratio | 0.2 |
| $\alpha$ | Biot's coefficient in the matrix | 0.8 |
| $\phi$ | Material porosity | 0.01 |
| $c_p$ | Fluid compressibility | $4.0 \times 10^{-10}\mathrm{Pa^{-1}}$ |
| $\mathcal{K}$ | Matrix permeability | $5.0 \times 10^{-20}\mathrm{m^2}$ |
| $\mu$ | Viscosity | $1.0 \times 10^{-4}$ Pa · s |
| $\mu_s$ | Friction coefficient | 0.5 |
| $\psi$ | Dilation angle | 1° |

The effects of the time step $(\Delta t)$, size of the microscale domain $(l)$, and mesh size $(\Delta H)$ on SIFs, displacement, and pressure are considered. Two levels of the microscale domain size are considered, i.e., $l = 0.5$ m or 1.0 m. The resolutions of the microscale and macroscale domains are the same, i.e., $\varepsilon_{\mathrm{m}} = 1.0$, and we consider two different levels: $\Delta H = \Delta h = 0.01$ m or $0.02$ m. Three levels of time steps are used: $\Delta t =$ 0.5 h, 1.0 h, or 1.5 h. Because of the lack of experimental data and exact solutions, the results from a computational setup are chosen as the reference. In the reference setup, the microscale domain coincides with the macroscale domain, i.e., $l = L$, and they are similar in resolution, i.e., $\Delta H =$ 0.01 m, $\varepsilon_{\mathrm{m}} = 1.0$. We also use a small time step, i.e., $\Delta t =$ 0.5 h, in this setup.

The SIFs at tip A, pressure at point $\mathrm{C} = (1.5, 1.5)$, and displacement at point $\mathrm{D} = (2, 2)$ obtained by different microscale domain sizes, mesh sizes, and time steps are shown in Fig. 6 and Fig. 7. There are no significant differences in the solutions obtained by using a small microscale domain compared to using a larger domain for a given mesh size $\Delta H$. In addition, the resolution of meshes and time steps have little effect on the solution. These agreements indicate that the calculation of SIFs is stable for all of the time step, microscale domain size, and mesh size considerations.

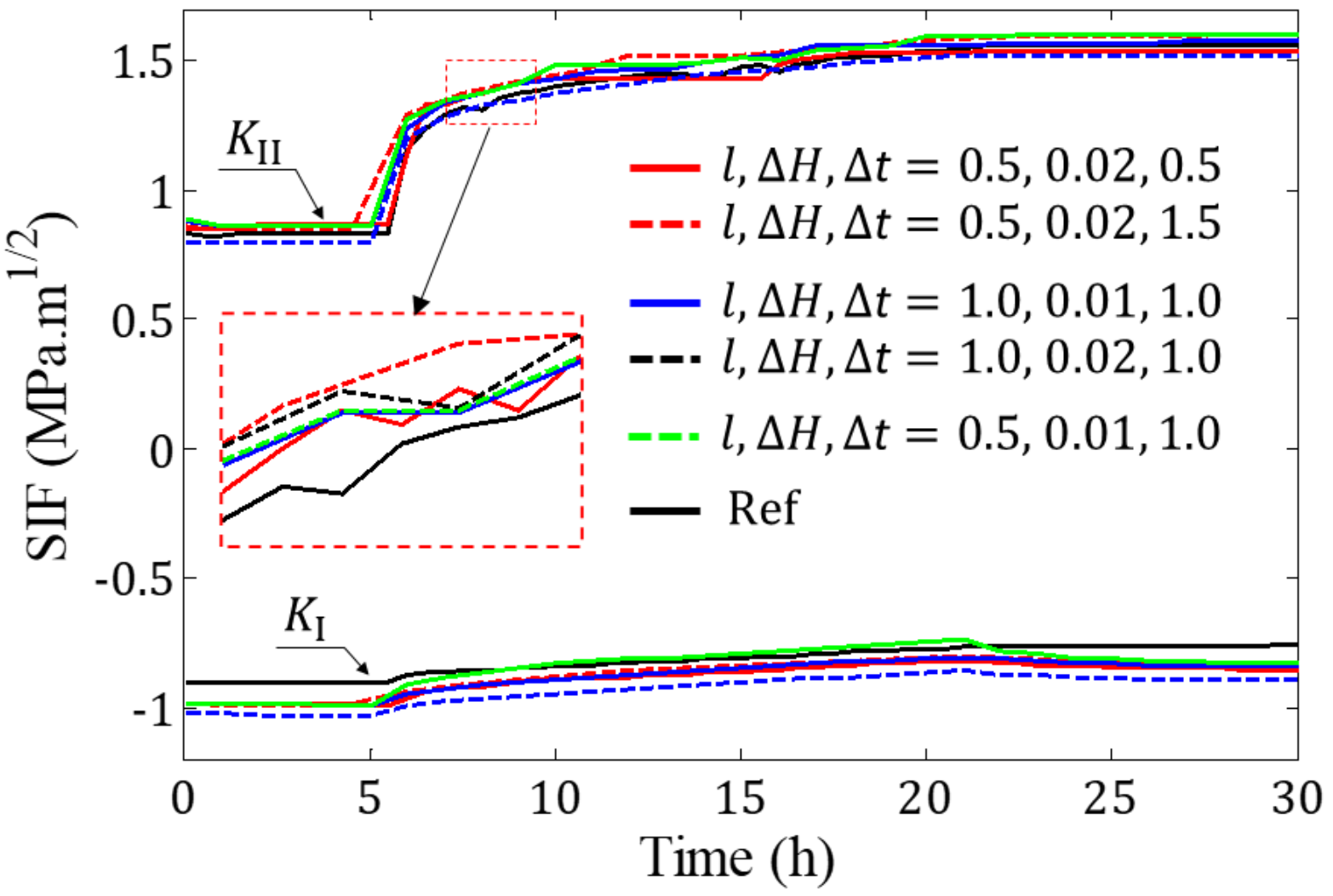

Fig. 6. Stress intensity factors at tip $A$ obtained by different microdomain sizes $l$ (m), mesh sizes $\Delta H$ (m), and time steps $\Delta t$ (h).

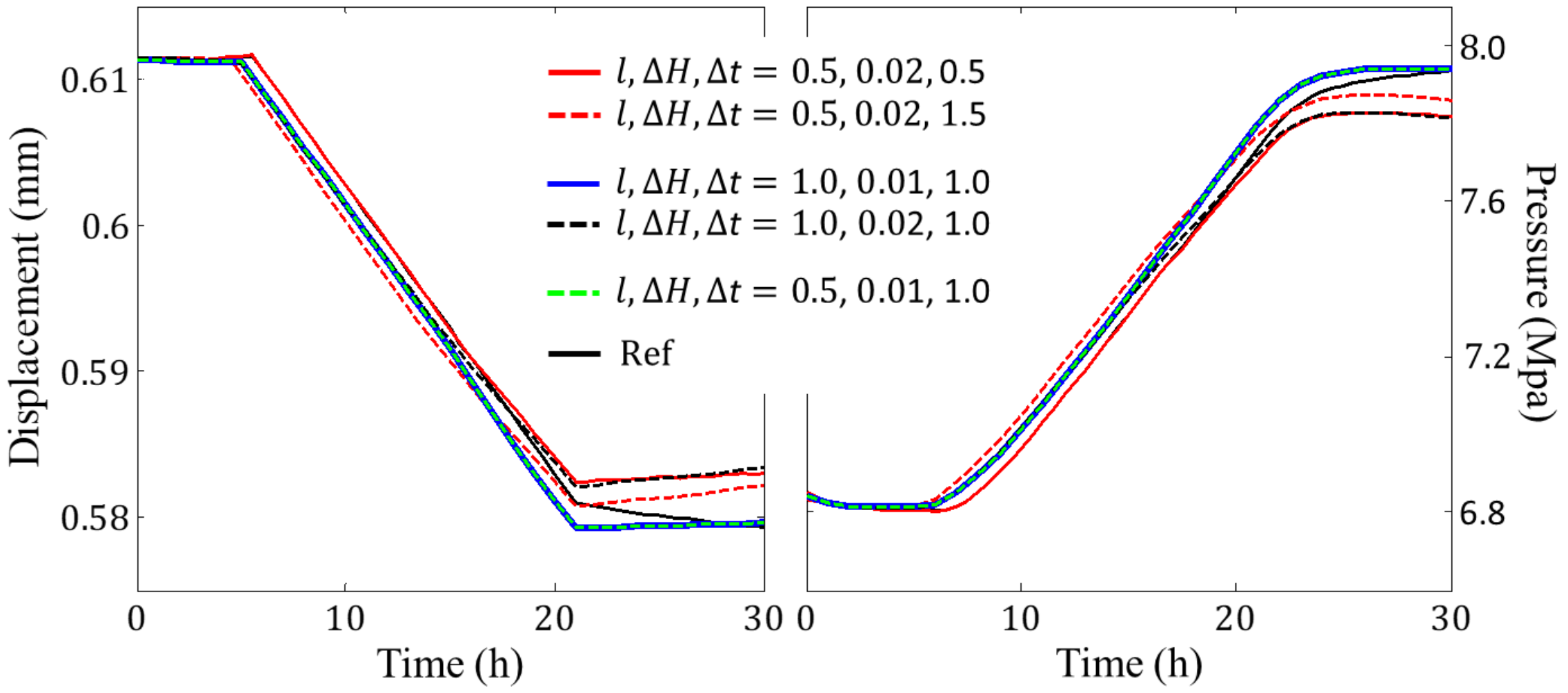

Fig. 7. Displacement at point $D$ and pressure at point $C$ obtained by different microdomain sizes $l$ (m), mesh sizes $\Delta H$ (m), and time steps $\Delta t$ (h).

The injection at a low rate gradually builds up pressure, causing slip of pre-existing fracture and shear failure instead of tensile failure at the fracture tip. As shown in Fig. 6, the injection increases $K_{\mathrm{II}}$ from $0.79\ \mathrm{MPa \cdot m^{1/2}}$ to $1.47\ \mathrm{MPa \cdot m^{1/2}}$, while there is almost no effect on $K_{\mathrm{I}}$. This is the result of the gradual reduction of the contact traction at the pre-existing fracture during the injection. In addition, the injection also increases pore pressure and resists deformation of the domain caused by compression. As shown in Fig. 7, the pore pressure at point C and displacement at point D are 6.8 MPa and 0.63 mm, respectively, at the beginning of the injection. After 15 h of injection, pressure increases to 8.0 MPa and displacement decreases to 0.60 mm. Termination of injection keeps the pore pressure and displacement stable.

### 5.2 Wing-crack propagation caused by fluid injection

Next, we consider further aspects of the model and solution strategy by studying the initiation and propagation of wing cracks from the ends of a pre-existing fracture caused by gradual pressure build-up by fluid injection at a low rate. The problem geometry, boundary conditions, and material parameters are the same as in the previous example. The fracture toughness is $K_{\mathrm{IC}} = 0.7\ \mathrm{MPa \cdot m^{1/2}}$. Water is injected at a constant rate of $Q_0 = 5 \times 10^{-9}\ \mathrm{m^2/s}$ into the pre-existing fracture during the simulation. The simulation is stopped when the wing crack propagating from tip A reaches a length of 0.25 m. To facilitate comparison between different sets of simulations, a reference case is computed with the same simulation parameters as used for the reference in Section 5.1, where we additionally use a threshold of $\varepsilon_{\mathrm{p}} = 0.5$ for the macroscale propagation of the fracture.

### 5.2.1 Simulation study of coupled physics

We start by illustrating the capacity of the present model to capture the complex coupled physics involved in fracture deformation and propagation based on a study of results from the reference case. Fig. 8 and Fig. 9 show the aperture expansion, shear slip, and contact traction of the fracture during the simulation. The fracture slowly extends in the first forty hours, but then it suddenly increases more quickly. Due to injection at a relatively low rate and low permeability of the matrix, pressure takes time to build up in the porous media domain. This process gradually reduces the fracture contact traction and causes small slips on the fracture surface, but the fracture is still in contact. The small slip causes slight growth along the pre-existing fracture. After approximately forty hours of injection, the fluid pressure is sufficiently elevated in the domain to decrease the contact traction and induce slip in larger regions of the pre-existing fracture, leading to the further propagation of wing cracks. Finally, the contact traction goes to zero along the initial fracture, leading to its complete opening and rapid propagation.

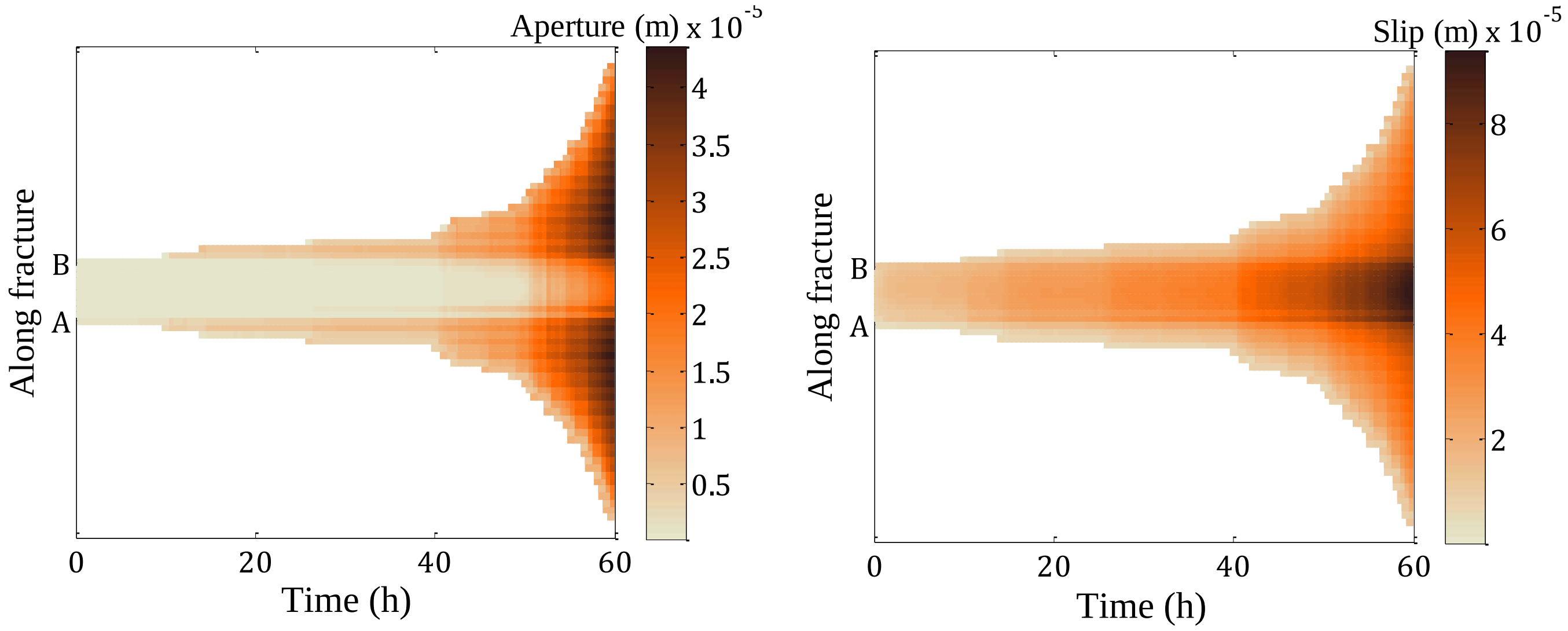

Fig. 8. Aperture expansion and shear slip at the fracture during simulation.

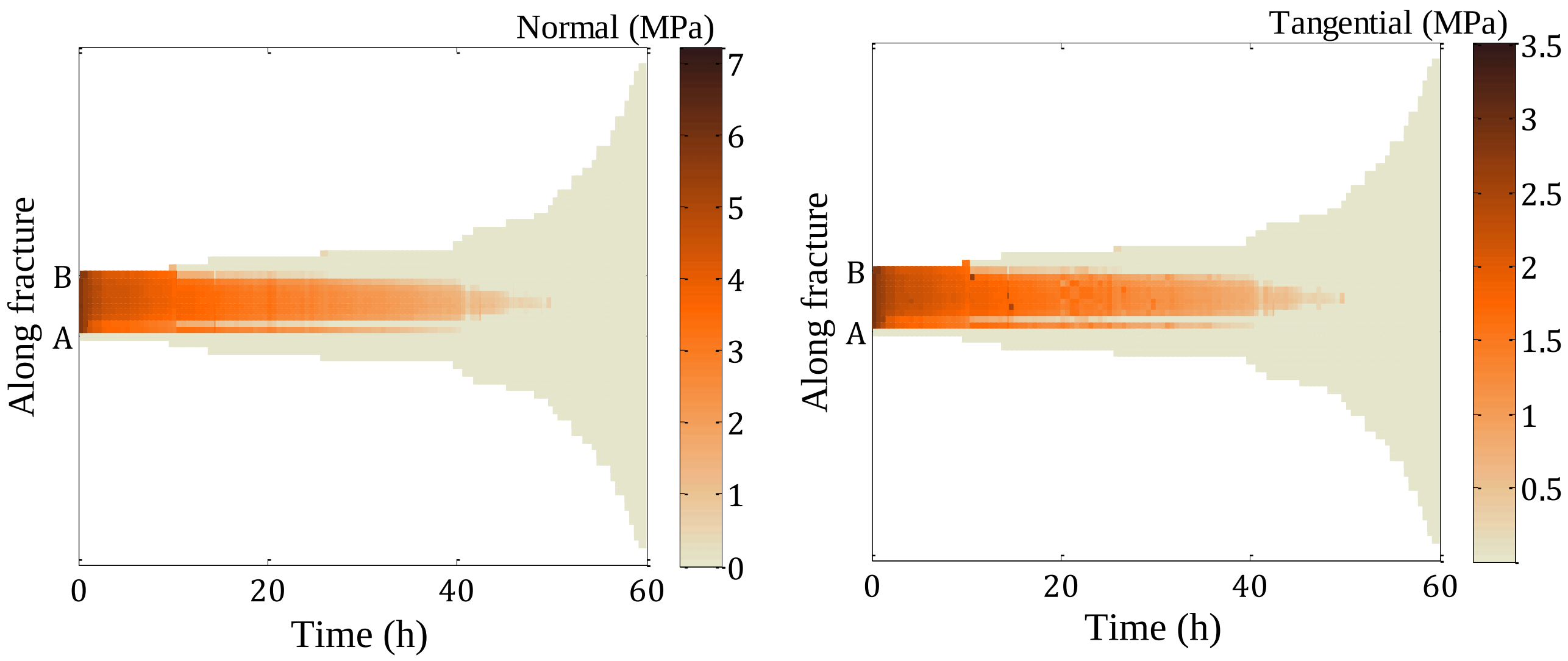

Fig. 9. Normal and tangential tractions at the fracture during simulation.

The fracture geometry and surrounding pore pressure for the simulations are shown in Fig. 10. First, the slip of fracture faces triggers wing cracks to initiate from the tips. Then, fractures slowly propagate toward the maximum horizontal stress direction (the *x*-direction) during the first forty hours. High fluid pressure mainly occurs along the pre-existing fracture due to the low permeability of the surrounding porous medium. After that, the fracture extension makes fluid pressure propagate further in the domain while injection continues to elevate the pressure. This process reduces contact traction and causes the fracture to grow farther.

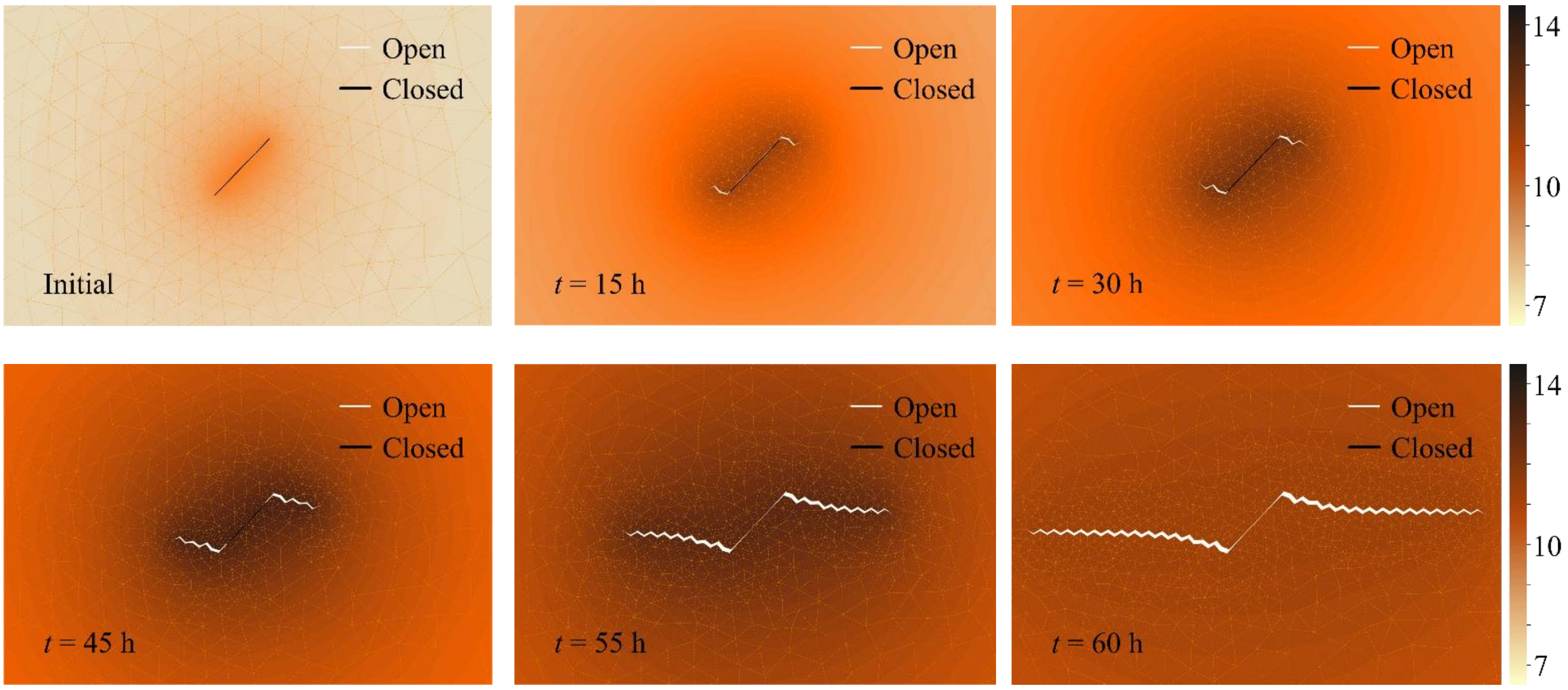

Fig. 10. Fracture propagation and pressure evolution in a 2D porous media during the fluid injection into the pre-existing fracture. Solid white/black lines denote fractures. The color bar represents pore pressure (MPa).

In the first ten hours, pressure and slip in the injected fracture increase gradually, not enough to trigger the fracture to propagate. So, the macroscale computational grid and solution are preserved. After that, wing cracks emerge and require the grid adjustment, followed by the update of the solution.

### 5.2.2 Effect of macroscale grid resolution

To probe the robustness of the multiscale simulation approach, we first investigate the effect of macroscale grid resolution on predicting the speed of fracture growth and fracture paths. In this example, the microscale grid coincides with the macroscale grid, i.e., $l = L$, and $\varepsilon_{\mathrm{m}} = 1.0$, while the time step is $\Delta t = 1.0$ h. A comparison of predictions obtained by different resolutions is shown in Fig. 11, in which black lines represent a prediction based on the computational reference.

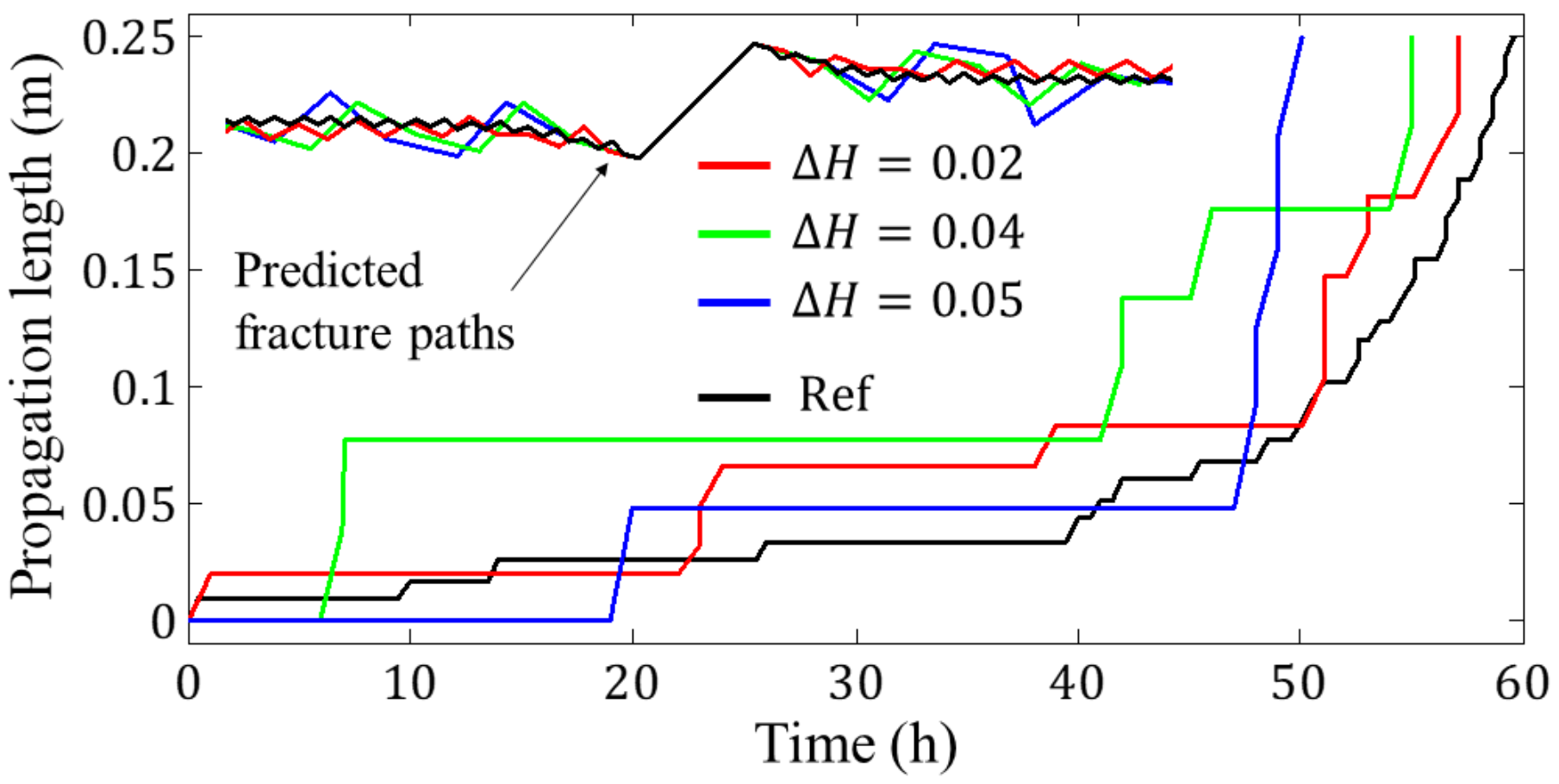

Fig. 11. The propagation length (m) from tip A and predicted fracture paths obtained by different resolutions $\Delta H$ (m) for $l = L$, $\varepsilon_{\mathrm{p}} = 1.0$, and $\Delta t = 1.0$ h

For all grid resolutions, the wing cracks propagate in the direction normal to the least principal stress, although the fracture paths can be seen to meander, particularly for the coarser grids. The propagation speed is initially stable with small increments in fracture size, followed by accelerated propagation starting at 40-50 h for the different grid resolutions. The results are in relatively good agreement in the first period, although the timing of the propagation events varies between the grid resolutions. The results differ more in the acceleration period, with the coarsest resolution ($\Delta H = \Delta h = 0.05$ m) showing almost brutal fracturing compared with the gradual although accelerating speed for the solutions obtained on the more refined grids. This is not unexpected since fast propagation is hard to capture, particularly for coarse grid resolutions. The example thus illustrates the balance between accuracy and computational cost and underlines the need to adapt and refine the macroscale grid in the vicinity of a propagating fracture tip.

### 5.2.3 Effect of microscale domain size and grid resolution

Next, we consider the impact of seeking computational savings in the microscale problem by assigning a smaller microscale domain size $l$ and different resolutions in the microscale grids. We fix the time step to $\Delta t = 1.0$ h and set the resolution of the macroscale domain and propagation threshold to $\Delta H = 0.02$ m and $\varepsilon_{\mathrm{p}} = 0.5$, respectively. The failure criterion is evaluated from a solution to the microscale problem on a domain surrounding the fracture tip.

The effect of the microscale domain size, $l$, and the resolution, $\Delta h = \varepsilon_{\mathrm{m}} \Delta H$, of this microscale domain on the propagation speed and fracture path is shown in Fig. 12. Again, the calculated fracture paths meander for the coarser grid, but this effect abates with the refinement of the microscale grid. In terms of propagation speed, the simulations again show a transition from stable to accelerating propagation. Except for $\varepsilon_{\mathrm{m}}$ in the period from approximately 8 h to 20 h after the start of the injection, the propagation speed is always larger for the smaller microscale domains. As the boundary conditions for the microscale problem are fixed by the macroscale state, the smaller domains must absorb the energy from fracture sliding in a smaller rock volume, increasing the stresses in the vicinity of the tip. The severity of this effect should, to a large degree, be independent of the size of the macroscale domain, and thus, using only somewhat larger microscale domains should also be feasible for larger problems.

The impact of varying the microscale mesh size is less clear. The results obtained on the two coarser microscale grids, $\varepsilon_{\mathrm{m}} = 0.5$ and $\varepsilon_{\mathrm{m}} = 1.0$, are in broad agreement in the period of stable propagation but exhibit notable differences when transitioning to an accelerating regime. The results from the finest microscale grid do not show a period of stable propagation on the macroscale grid but instead enter the period of acceleration directly, at approximately the same time as the other simulations. This disagrees with the other results, notably the observation in Fig. *6* that the critical threshold for $K_{\mathrm{II}}$ is crossed approximately 7 h after the start of injection. A possible explanation is that the relatively large difference in mesh size between the microscale and macroscale problems for this value of $\varepsilon_{\mathrm{m}}$ makes the fracture propagate on the microscale without this effect being captured on the macroscale domain.

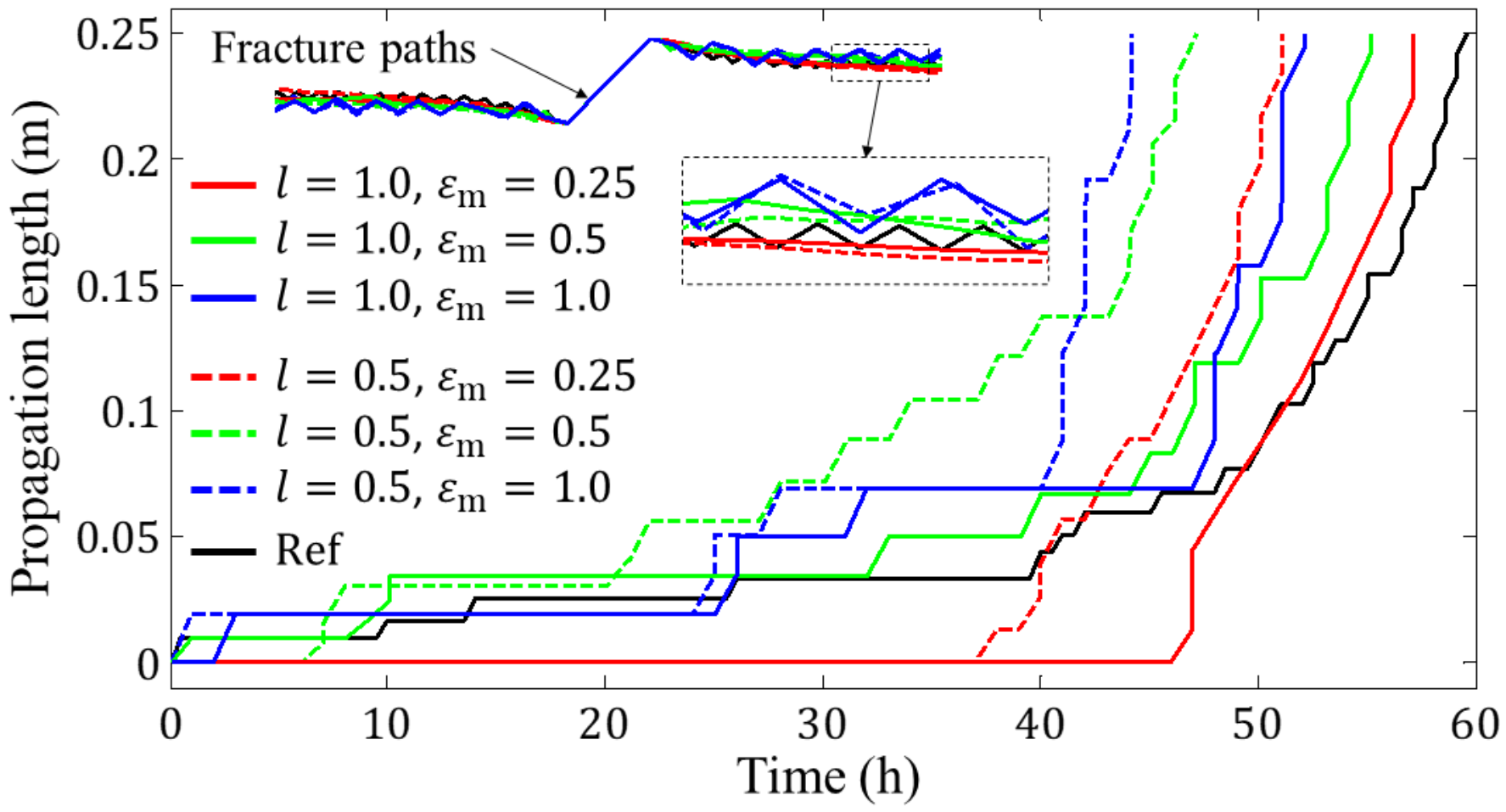

Fig. 12. The propagation length from tip A and fracture paths obtained by different sizes of the microscale domain, $l$ (m), and resolutions of the microscale domain, $\varepsilon_m$, for $\Delta H = 0.02\ m,\ \Delta t = 1.0\ h,$ and $\varepsilon_p = 0.5.$

### 5.2.4 Effect of discretization coupling parameters: Time-step size and macroscale propagation threshold

The effects of time step and threshold to extend fracture in the macroscale domain are investigated in this example. For the time step study, as shown in Fig. 13 (a), we fix the size of the microscale domain and mesh resolutions, i.e., $l = 1.0$ m and $\Delta H = \Delta h = 0.01$ m, and consider three levels of the time step, i.e., $\Delta t = 0.5$ h, $1.0$ h, or 1.5 h. The method can be seen to be stable under this variation with a similar speed of fracture propagation in all three cases. The difference from the reference case can be attributed to the smaller size of the microscale domain.

Besides, we also consider how the macroscale propagation threshold, $\varepsilon_{\mathrm{p}}$, affects the propagation speed, as shown in Fig. 13 (b). In this case, we set $l = 1.0$ m, $\Delta H = 0.02$ m, $\Delta t = 1.0$ h, and $\varepsilon_{\mathrm{m}} = 1/3$ and assign three different values for $\varepsilon_{\mathrm{p}}$, namely, 1/3, $2/3$, or 1. For the two highest values of $\varepsilon_{\mathrm{p}}$, i.e., $2/3$ and $1$, there is no period of stable propagation on the macroscale but rather an abrupt transition to the accelerating regime. This is like the results reported in Fig. *12*, which also had a small value of $\varepsilon_{\mathrm{m}} = 1/4$. Using a small value for $\varepsilon_{\mathrm{p}}$ also compensates for this effect.

The findings of these experiments are summarized as follows. First, the direction of propagation is mainly controlled by boundary conditions. However, the computed fracture path tends to wiggle unless a relatively fine grid is applied on the microscale domain. Such refined grids may again lead to delayed propagation on the macroscale grid unless the microscale and macroscale problems are tightly coupled through the parameter $\varepsilon_{\mathrm{p}}$. Second, the size of the microscale domain can be reduced to lower the computational cost. However, a microdomain that is too small can overestimate the propagation speed. Finally, the time step size had little impact on the results for the cases we considered.

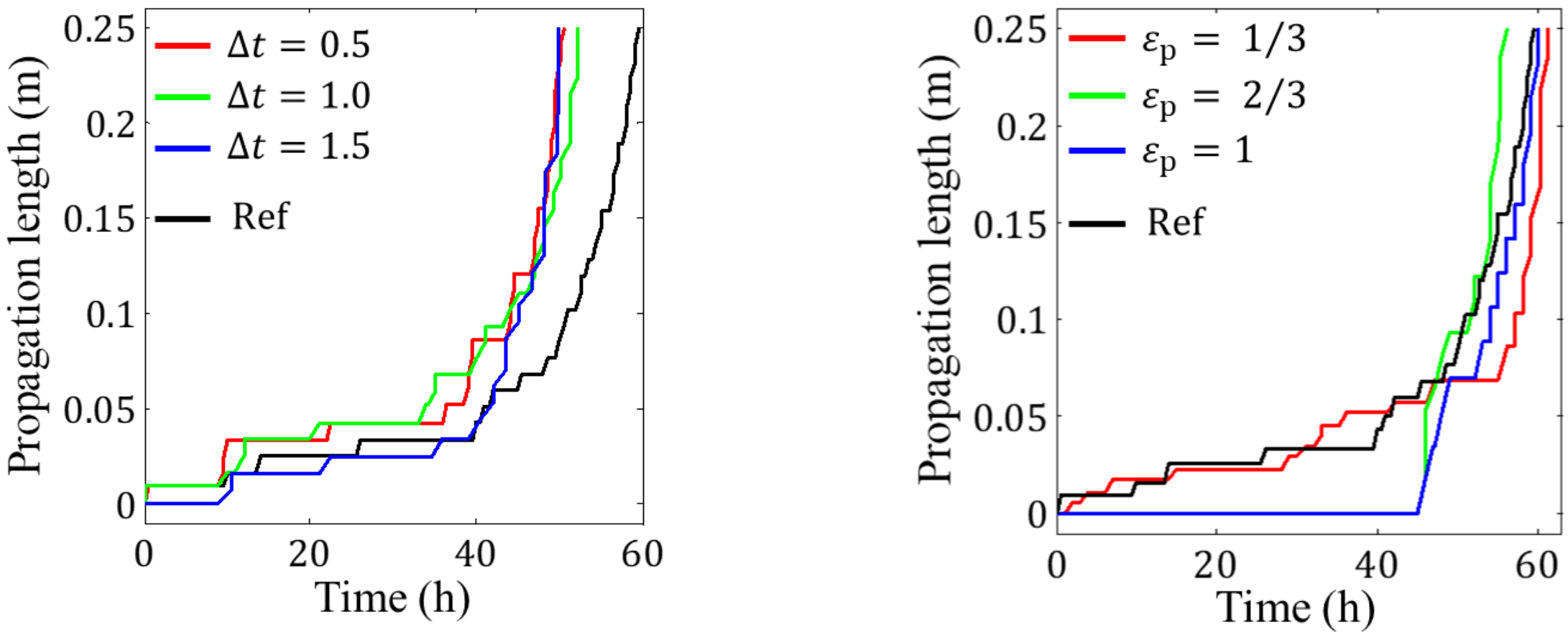

a) Effect of time step $\Delta t$ (h) for $l = 1.0$ m, $\Delta H = 0.01$ m, $\varepsilon_{\mathrm{m}} = 1.0$, and $\varepsilon_{\mathrm{p}} = 1.0$

b) Effect of propagation threshold $\varepsilon_{\mathrm{p}}$ for $l = 1.0$ m, $\Delta H = 0.02$ m, $\varepsilon_{\mathrm{m}} = 1/3$, and $\Delta t = 1.0$ h

Fig. 13. The propagation length from tip A obtained by different resolutions, microscale domain sizes, time steps, and propagation thresholds.

### 5.3 Extension: Propagation of multiple fractures under fluid injection

Finally, to show the power of the proposed approach, hydromechanical processes interacting with the deformation and propagation of three pre-existing fractures in porous media are studied. The geometry and boundary conditions are illustrated in Fig. 14, in which $L_x = L_y = 2$ m. There are three fractures with the same initial aperture, $a_i^0 = 1$ mm, pre-existing in the domain. Fracture 1 is defined by endpoints $\mathrm{A} = (0.751, 1.208)$ and $\mathrm{B} = (0.849, 1.192)$, fracture 2 by endpoints $\mathrm{C} = (0.965, 0.965)$ and $\mathrm{D} = (1.035, 1.035)$, and fracture 3 by endpoints $\mathrm{E} =$ $(1.144, 0.780)$ and F $= (1.256, 0.831)$. Water is injected at a constant rate of $Q_0 = 5 \times 10^{-9}\ \mathrm{m}^2/\mathrm{s}$ into fracture 2 during the simulation. The material parameters are the same as for the example in Section 5.2. The simulation is implemented based on the multiscale model, in which $l = 0.5$ m, $\Delta H = 0.02$ m, $\varepsilon_{\mathrm{m}} = 0.5$, $\varepsilon_{\mathrm{p}} = 0.5$, and $\Delta t = 1$ h.

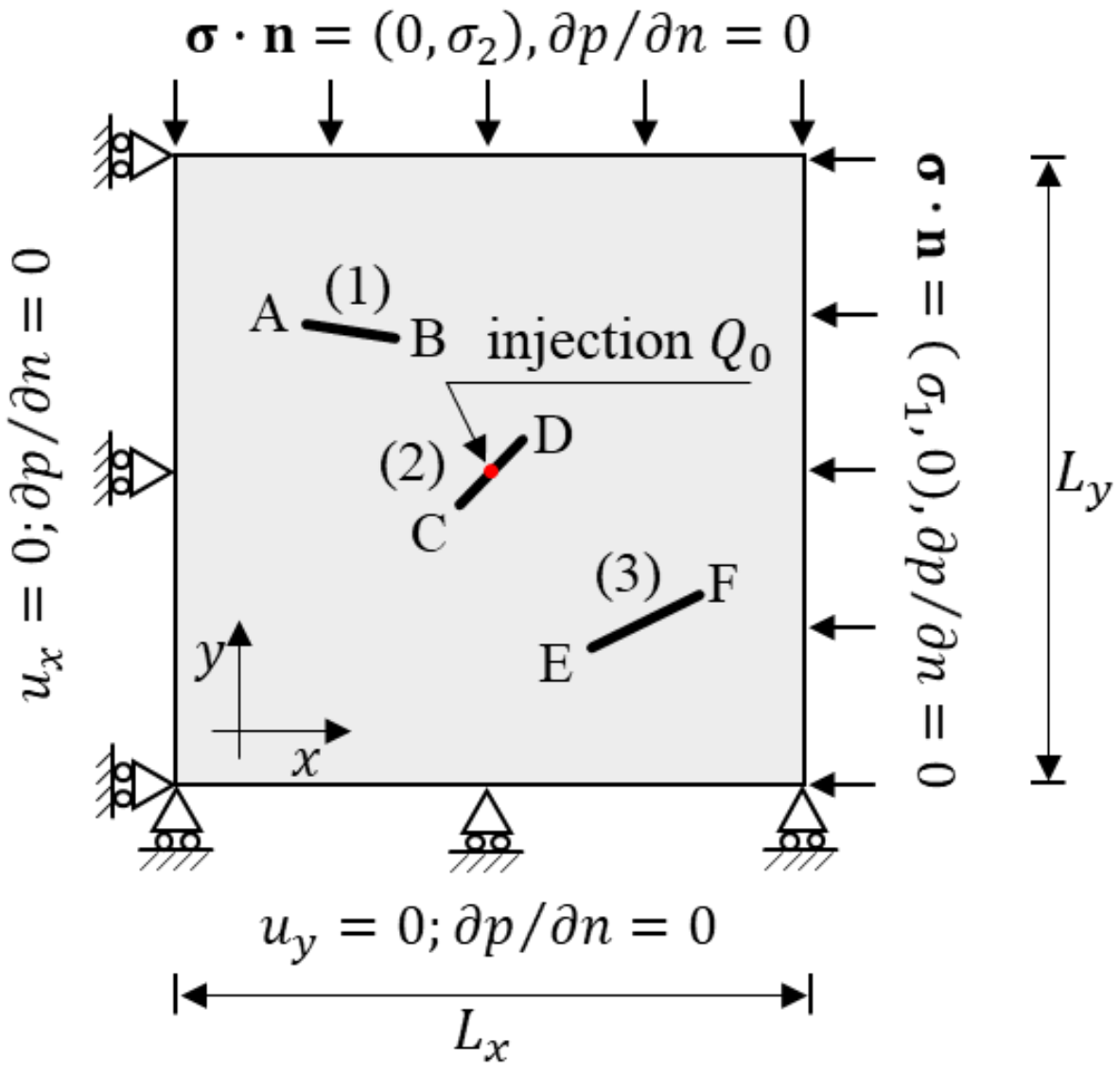

Fig. 14. Model of three fractures in a porous media subject to boundary conditions.

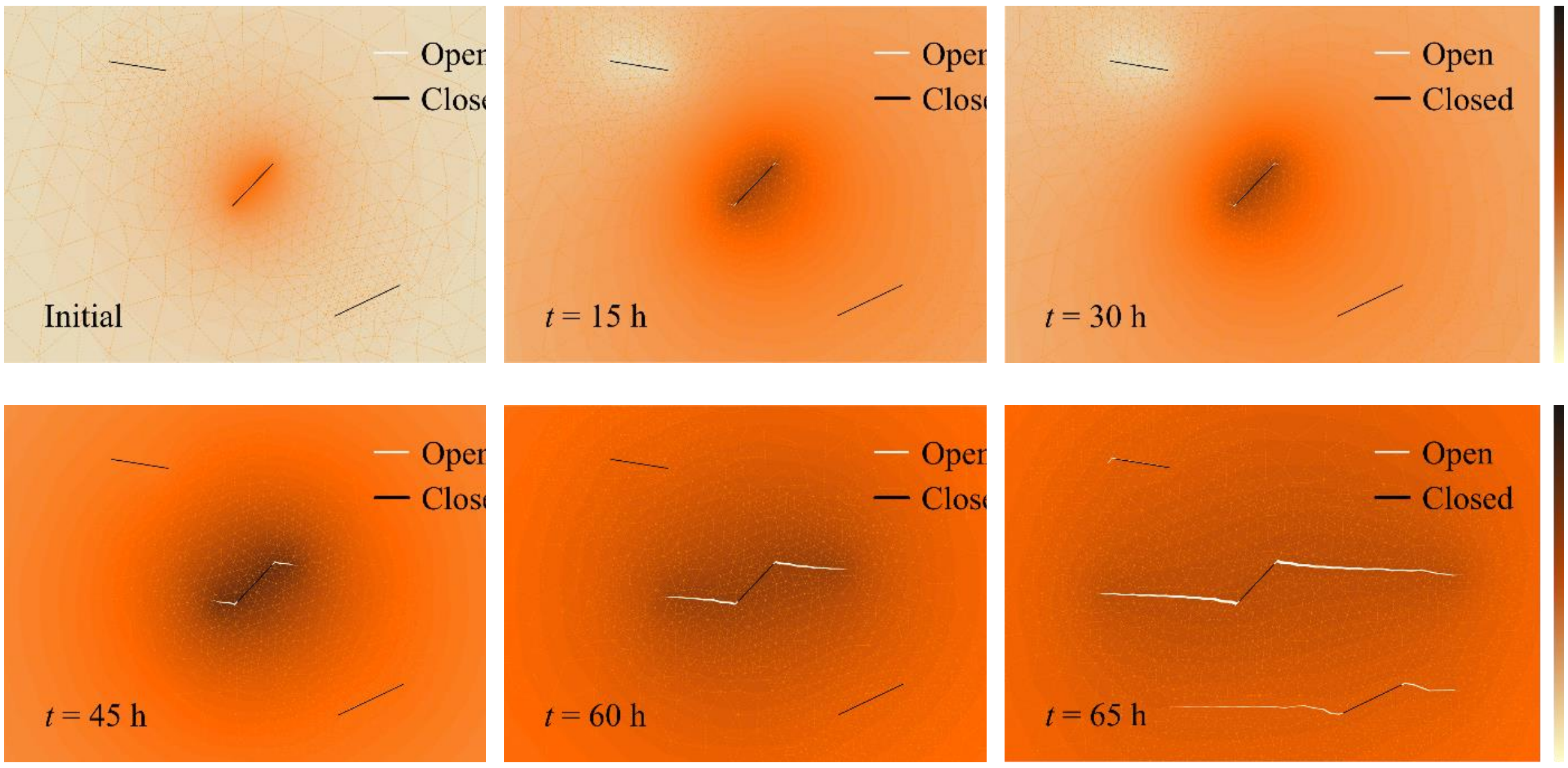

Fig. 15. Propagation of multiple fractures and evolution of pressure in a 2D porous media during fluid injection into one pre-existing fracture. Solid white/black lines denote fractures. The color bar represents pore pressure (MPa).

The evaluation of the fracture geometry and the pore pressure is shown in Fig. 15. As in the above example, wing cracks mainly initiate and propagate from the fracture where fluid is injected. Small wing cracks are also observed from the other pre-existing fractures and are caused mainly by mechanical effects, i.e., the deformation of the domain. The fractures where fluid is injected propagate as a consequence of the hydromechanical stresses induced by the fluid injection, and fluid infiltrates farther into the domain as the fracture grows. These processes take place at the same time, leading the fracture to grow increasingly faster. In addition to the fracture growth, the injection also stimulates the nearby fracture to propagate.

## 6. Conclusion

This work presented a mathematical model and a numerical solution for coupling fluid flow, matrix deformation, fracture slip, and fracture propagation in porous media due to fluid injection. The governing mathematical model is based on Biot's model, with the deformation of existing fractures represented by contact mechanics. The maximum tangential stress criterion is combined with Paris' law to govern the fracture growth processes. A multiscale simulation was presented to reduce computational cost and ensure accuracy. The numerical approach employs a novel combination of finite volume methods for the poroelastic deformation of existing fractures with a finite element approach for the fracture propagation process.

The verifications in this paper show that the proposed approach is stable with different time steps, macroscale grid sizes, and microscale grid sizes. This approach is capable of simulating complex problems, such as the simultaneous propagation of multiple fractures combined with the slip and dilation of fractures in contact and tensile opening. Hydraulically and mechanically interacting fractures are handled naturally. Therefore, the current model has potential in the simulation of mixed-mechanism hydraulic stimulation of fractured reservoirs, in which both fracture shearing and corresponding wing-crack propagation lead to an increase in permeability.

## Declaration of Competing Interest

The authors declare that they have no known competing financial interests or personal relationships that could have appeared to influence the work reported in this paper.

## Acknowledgments

This project has received funding from the European Research Council (ERC) under the European Union's Horizon 2020 research and innovation programme (grant agreement No 101002507).